\documentclass[reqno,twoside,final]{amsart}
\usepackage{hyperref}
\usepackage{amssymb,amsmath}

\newtheorem{Conj}{Conjecture}
\newtheorem{Thm}{Theorem}
\newtheorem{Cor}{Corollary}
\newtheorem{Lem}{Lemma}

\theoremstyle{definition}
\newtheorem{Def}{Definition}
\theoremstyle{remark}

\numberwithin{equation}{section}

\newcommand{\thmref}[1]{Theorem~\ref{#1}}
\newcommand{\secref}[1]{\paragraph\ref{#1}}
\newcommand{\lemref}[1]{Lemma~\ref{#1}}

\newcommand{\defref}[1]{Definition~\ref{#1}}
\newcommand{\corref}[1]{Corollary~\ref{#1}}

\def\e{\varepsilon}

\def\R{{\mathbb R}}

\def\G{\varGamma}
\def\F{F_{\varepsilon}}
\def\K{K_{\varepsilon}}
\let\paragraph=\S
\let\hataccent=\^
\def\S{\mathbb S}

\def\({\left(}
\def\){\right)}
\def\[{\left[}
\def\]{\right]}

\def\const{\operatorname{const}}
\def\dist{\operatorname{dist}}
\def\maxdist{\operatorname{maxdist}}

\def\codim{\operatorname{codim}}

\def\sign{\operatorname{sign}}
\def\k{\operatorname{\Bbbk\kern1pt}}

\def\ud{\mathrm{d}}

\def\:{\colon}
\def\~#1{\widetilde{#1}}
\def\^#1{\widehat{#1}}
\def\=#1{\check{#1}}
\let\le=\leqslant
\let\ge=\geqslant

\def\<{\left<}
\def\>{\right>}
\begin{document}
\title[About affine version of Arnold Conjecture]{On affine hypersurfaces
with everywhere nondegenerate Second Quadratic Form}
\author{A. Khovanskii, D. Novikov}
\date{\today}
\address{Department of Mathematics, Toronto University, Toronto, Canada}
 \email{dmitry@math.toronto.edu}
\address{Department of Mathematics, Toronto University, Toronto, Canada}
 \email{askold@math.toronto.edu}
\thanks{Khovanskii's work is partially supported by
Canadian Grant {\rm N~0GP0156833}. Novikov's work was supported by
the Killam grant of P. Milman and by James S. McDonnell
Foundation.}
\begin{abstract}
Consider a closed connected hypersurface in $\mathbb{R}^n$ with
constant signature $(k,l)$ of the second quadratic form, and
approaching a quadratic cone at infinity.  This hypersurface
divides $\mathbb{R}^n$ into two pieces. We prove that one of them
contains a $k$-dimensional subspace, and another contains a
$l$-dimensional subspace, thus proving an affine version of Arnold
hypothesis. We construct an example of a surface of negative
curvature in $\mathbb{R}^3$ with slightly different asymptotical
behavior for which the previous claim is wrong.
\end{abstract} \maketitle
\tableofcontents

\section{Introduction}
In this paper we prove three results connected to an Arnold
hypothesis formulated in \cite{Arnold88}. Recall first this
hypothesis.

\begin{Def}\label{(k,l)-hyperbolic}
A smooth hypersurface in $\mathbb{RP}^{n+1}$ is called
$(k,l)$-hyperbolic if its second quadratic form has signature
$(k,l)$. In other words, in some affine system of coordinates
$(x,y,z)$ its local  equation is
$z=\sum_{i=1}^{k}x_{i}^2-\sum_{j=1}^{l}y_i^2 +\text{higher order
terms}$.
\end{Def}

One class of examples of such hypersurfaces is the class of
quadrics. Namely, let $Q(x,y)$ be  a non-degenerate symmetric
bilinear form in  ${\mathbb R}^{n+1}$ of signature $(k+1,l+1)$
(i.e. its restriction to some $(k+1)$-dimensional subspace $L_+$ is
positively defined, and its restriction to some
$(l+1)$-dimensional subspace $L_-$ is negatively defined, and
$k+l=n-1$). Then a hypersurface $S_Q$ in ${\mathbb R}P^n$ given by
equation $Q(x,x)=0$ is smooth and $(k,l)$-hyperbolic.

The hypersurface $S_Q$ has the following remarkable property: the
domain bounded by $S_Q$ (i.e. the domain $\{Q(x,x)\le 0\}$)
contains a $l$-dimensional projective subspace $P(L_+)$ -- a
projectivization of $L_+$ , and its complement contains a
$k$-dimensional projective subspace $P(L_-)$ -- a projectivization
of $L_-$ .

Arnold's hypothesis claims existence of such subspaces for any
$(k,l)$-hyperbolic hypersurface.
\begin{Conj}[Arnold Conjecture]
\noindent 1. For any domain $U\subset{\mathbb R}P^n$  bounded by a
connected smooth $(k,l)$-hyperbolic hypersurface $B$ there exist a
projective subspace $L_-$ of dimension $k$ not intersecting $U$
and a projective subspace $L_+$ of dimension $l$ contained in $U$.

\noindent 2. Any projective line joining a point of $L_+$ and a
point  $L_-$ intersects $B$ at exactly 2 points.
\end{Conj}

Apart from the case of quadrics discussed above, another fact
justifying this conjecture is the following well-known fact.
Consider a locally convex connected surfaces $S$ in
$\mathbb{R}P^n$. Then $S$ is a boundary of  some  domain in
$\mathbb{R}P^n$ and this domain doesn't intersects some
hyperplane, see \cite{Arnold88}. This case corresponds to the case
of $k=0$  or $l=0$ of the Arnold Conjecture.

We prove in \cite{NKh1} the first case $k=l=1$ of the Arnold
Conjecture in some additional assumptions. Namely, we consider a
smaller then class of $(1,1)$-hyperbolic surfaces class of {\em
projective} $L$-convex-concave subsets of ${\mathbb R}P^3$. This
class is a close relative of the class convex-concave sets
considered in \secref{ssec:affine c-c sets}. We prove that any
$L$-convex-concave subset contains a line inside.

In this paper we deal with affine version of Arnold's hypothesis.
Namely, we consider $(k-1,l)$-hyperbolic hypersurfaces in
${\mathbb R}^{k+l}$ with some prescribed asymptotic behavior at
infinity.  Our results can be roughly summarized as follows: if
the asymptotic at infinity forces the closure in ${\mathbb
R}P^{k+l}$ of the hypersurface to be $(k-1,l)$-hyperbolic, then
the domain bounded by the hypersurface contains a line. And if the
closure is not $(k-1,l)$-hyperbolic, then one can construct a
domain bounded by such hypersurface and not containing  any line
inside.

Here is more exact description of the results. Consider
$(k-1,l)$-hyperbolic hypersurfaces $M$ in ${\mathbb R}^n$,
$k+l=n$. We say that $M$ approaches a surface $L$ at infinity if
$M$ and $L$ are arbitrarily $C^2$-close outside a big enough ball
(see \secref{ssec:Rolle definitions} for an exact
$\epsilon,\delta$-definition). For example, the quadric
$\{\sum_{i=1}^k x_i^2-\sum_{j=1}^{l} x_{k+j}^2\}=A$ approaches the
cone $K=\{\sum_{i=1}^k x_i^2=\sum_{j=1}^{l} x_{k+j}^2\}\subset
{\mathbb R}^n$ at infinity for any $A$.

We prove the following theorem:

\begin{Thm}\label{thmmain:1 line exists}
The first claim of the Arnold hypothesis is true for any
$(k-1,l)$-hyperbolic hypersurface $M$ approaching a quadratic cone
$K=\{\sum_{i=1}^k x_i^2=\sum_{j=1}^{l} x_{k+j}^2\}\subset {\mathbb
R}^n$ at infinity.
\end{Thm}

We also prove the following result, strengthening the second part
of the Arnold Conjecture.

\begin{Thm}\label{thmmain:2 projection is in}
Any $(1,1)$-hyperbolic surface $M$ in ${\mathbb R}^3$ approaching
a quadratic cone $K=\{x^2+y^2=z^2\}$ intersects any ray from the
origin in at most one point. In other words, the projection $M\to
\mathbb{S}^2=\{x^2+y^2+z^2=1\}$ is embedding.
\end{Thm}

Note that the surface described in \thmref{thmmain:1 line exists}
has a $C^1$-closure in $\mathbb{R}P^n$. Consider the simplest case
$k=l=1$, and denote $K_-=K\cap \{z\le 0\}$ and $K_+=K\cap \{Z\ge
0\}$. Will the result remain true if we consider surface $B$
approaching union of translated $K_-$ and $K_+$? It turns out that
if the translates intersect, then the domain bounded by the
surface $B$ still contains a line, and if the intersection is
empty, then this is not necessarily true. Note the different
behavior of the projective closure $\overline{B}$ of the surface
in these two cases. In both cases the points of
$\overline{B}\setminus B$ are not smooth points of the closure.
However, in the first case the  $\overline{B}$ is
$(1,1)$-hyperbolic after arbitrarily small perturbation, and in
the second case the surface $\overline{B}$ is locally convex at
these points. This coincides with what the Arnold Conjecture
prescribes, further strengthening it.
\begin{Thm}\label{thmmain:3 example}
Let $K'=\{(x,y,z)\quad|\quad x^2+y^2=(|z|-1)^2, |z|\ge 1\}$ be a
union of non-intersecting translates of $K_-$ and $K_+$. There are
$(1,1)$-hyperbolic surfaces approaching $K'$ at infinity and not
containing any lines.
\end{Thm}

All these result can be considered in more general context of
existence and properties of a solution of some boundary problem. A
natural boundary problem is to find a compact smooth
$(1,1)$-hyperbolic surface which boundary is a given tuple of
non-intersecting closed smooth curves and which is tangent at this
boundary to the given set of planes tangent to the curves. The
\thmref{thmmain:1 line exists} follows from the fact that solution
of some boundary problem of this type cannot intersect an open
domain -- the interior of the cone $K$.

\section{Preliminaries: Quadrics} The hypersurface $M$ of the
\thmref{thmmain:1 line exists} is approaching at infinity a cone
given by a quadric. In this paragraph we collect some standard
facts about quadrics we will need later.

\begin{Thm}\label{thm:quadrics} Let
$Q(x)=\sum_{i=1}^{k}x_{i}^2-\sum_{j=1}^{l}x_{k+j}^2$, and let
$Q_{\e}=\{f=\e\}$ be its level hypersurfaces.
\begin{enumerate}
\item restriction of $Q$ to a tangent plane $T_xQ_{\e}\subset\mathbb{R}^{k+l}$
at point $x\in Q_(\e)$ is a polynomial of second degree having
signatures $(k, l-1)$, $(k-1,l-1)$ or $(k-1,l)$ according to
$\e<0$, $\e=0$ or $\e>0$ correspondingly;
\item $Q_{\e}$ is $(k, l-1)$-hyperbolic if $\e<0$ and $(k-1, l)$-hyperbolic if $\e>0$;
\item Projectivization of $Q_0$ is a $(k-1,l-1)$-hyperbolic
hypersurface in ${\mathbb R}P^{k+l-1}$.
\end{enumerate}
\end{Thm}
\begin{proof}

First we prove a general Lemma. To a nondegenerate bilinear form
$q$ on a linear space $L$ corresponds an isomorphism
$\widetilde{q}$ between $L$ and its dual $L^*$, defined by the
condition $\ell(x)=q(\widetilde{q}(\ell),x)$ for $\ell\in L^*$ and
all $x\in L$. Using this form one can define a bilinear form
$q^*(\ell_1,\ell_2)=q(\widetilde{q}(\ell_1),\widetilde{q}(\ell_2))$
on $L^*$.

\begin{Lem}\label{lem:bilinear form} Restriction
of $q$ on a hyperplane $H=\{\ell=0\}\subset L$ has signature
$(k,l-1)$, $(k-1,l-1)$, or $(k-1,l)$ if $q^*(\ell,\ell)<0$,
$q^*(\ell,\ell)=0$ or $q^*(\ell,\ell)>0$ correspondingly.
\end{Lem}
\begin{proof}
By definition $q(\widetilde{q}(\ell),x)=0$  is equivalent to $x\in
H$. So $q^*(\ell,\ell)=0$ means $y=\widetilde{q}(\ell)\in H$, i.e.
the restriction of $q$ on $H$ is degenerate. From the other hand,
$H$ should intersect subspaces where $q$ is positive/negative
definite by subspaces of dimensions $k-1$ and $l-1$ at least, so
the only possible signature is $(k-1,l-1)$. Vice versa, if for
some $y\in H$  and for all $x\in H$ we have $q(y,x)=0$, then
$\ell$ is proportional to $\widetilde{q}^{-1}(y)$ and so
$q^*(\ell,\ell)=q(y,y)=0$.

So suppose $q^*(\ell,\ell)\neq 0$. Then Gramm-Schmidt procedure
starting from $y=\widetilde{q}(\ell)$ is possible and results in a
basis which first vector is $y$ and all the rest  span a
hyperplane of vectors $q$-orthogonal to $y$, i.e. $H$. Since in
this basis $q$ is diagonal, we easily see that the signature of
the restriction of $q$  to $H$ is as required.
\end{proof}

To prove the \thmref{thm:quadrics} apply \lemref{lem:bilinear
form} to $Q(x)=q(x,x)$ in ${\mathbb R}^n$, where $n=k+l$. The
first claim of the \thmref{thm:quadrics} follows from the
\lemref{lem:bilinear form} and the fact that
$\widetilde{q}(dQ(x))=2x$, so
$q^*(dQ(x),dQ(x))=4q(x,x)=4\epsilon$.

To prove the second claim of the \thmref{thm:quadrics} change
coordinates in such a way that $x\in Q_(\e)$ is an origin and
$T_xQ_(\e)=\{x_n=0\}$. The restriction $Q_1$ of $Q$ to $T_xQ_(\e)$
is a quadratic polynomial,  without linear and free terms in these
coordinates. So $Q=\e+\lambda x_n+Q_1(x_1,...,x_{n-1})
+x_n\ell(x_1,...,x_{n-1})$, where $Q_1$ is a homogeneous
polynomial of degree 2 and $\ell(x_1,...,x_{n-1})$ is linear. An
easy computation shows that
$x_n=-\lambda^{-1}Q_1(x_1,...,x_{n-1})+...$ on $\{Q=\e\}$. The
second claim follows now from the first claim of the
\thmref{thm:quadrics}.

The third claim follows since the kernel of the bilinear form
$q_1$ on $T_xQ_(\e)$ corresponding to the quadratic form $Q_1$ is
exactly the kernel of the projection ${\mathbb R}^{n}\to{\mathbb
R}P^{n-1}$.

\end{proof}

 Denote by ${\mathbb S}^{k+l-1}$ the standard sphere
$\{\sum_{i=1}^{k}x_{i}^2+\sum_{j=1}^{l}x_{k+j}^2=1\}\subset
\mathbb{R}^{k+l}$. The standard scalar product $(x,x')=\sum_i
x_ix_i'$ in $L$ defines, as above, a scalar product in $L^*$ and
also an isomorphism of $L$ and $L^*$. For a smooth cooriented
hypersurface $M\subset{\mathbb R}^n$ the Gauss map $G:M\to
{\mathbb S}^{n-1}$ maps a point $x\in M$ to the vector normal to
$M$ at $x$. If $M=\{P=0\}$ and $dP\not =0$ on $M$, then the Gauss
map is a composition of the map $x\to \frac {dP}{\|dp\|}$ and the
isomorphism between ${\mathbb R}^n$ and $({\mathbb R}^n)^*$
provided by a standard scalar product. A classical computation
shows that the Jacobian of the Gauss map is exactly the Gaussian
curvature of the hypersurface.
\begin{Lem}\label{lem:quadric's Gauss} Let
$Q=\sum_{i=1}^{k}x_{i}^2-\sum_{j=1}^{l}x_{k+j}^2$ be as in
\thmref{thm:quadrics}. Then the pull-back of  $Q^*$ from
$({\mathbb R}^n)^*$ to ${\mathbb R}^n$ by the isomorphism given by
the standard scalar product coincides with $Q$. The Gauss map
provides diffeomorphisms between $\{Q=\e>0\}$ and ${\mathbb
S}^n\cap \{Q>0\}$, between $\{Q=\e<0\}$ and ${\mathbb S}^n\cap
\{Q<0\}$, and maps ${\mathbb S}^n\cap \{Q=0\}$ diffeomorphically
onto itself
\end{Lem}
This follows immediately from the explicit formulae for the Gauss
mapping  $G$ of the quadric ${\{Q=\e\}}$. Namely,
$$G(x_1,...,x_k,x_{k+1},...,x_{k+l})=\frac
1{\|x\|}(x_1,...,x_k,-x_{k+1},...,-x_{k+l}).$$ \qed

The topology of the domains, and therefore of level hypersurfaces
of $Q$ is quite simple.
\begin{Cor}\label{cor:quadric's topology}
$\{Q=\e>0\}$ is diffeomorphic to $\S^{k-1}\times B^l$,
$\{Q=\e<0\}$ is diffeomorphic to $\S^{l-1}\times B^k$ and
${\mathbb S}^n\cap\{Q=0\}$ is diffeomorphic to $\S^{k-1}\times
\S^{l-1}$
\end{Cor}
For example, ${\mathbb
S}^n\cap\{Q=0\}=\{\sum_{i=1}^{k}x_{i}^2=\sum_{j=1}^{l}x_{k+j}^2=\frac
1 2\}=\S^{k-1}\times \S^{l-1}$.

\section{Surface doesn't enter the half-cones}\label{sec:Rolle}
\subsection{Exact formulation of the result}\label{ssec:Rolle definitions}
In this section we prove the \thmref{thmmain:1 line exists}. We
start with a definition of an affine version of
$(k,l)$-hyperbolicity and state more precisely the asymptotic
conditions on the hypersurface $M$.
\begin{Def}
A smooth connected hypersurface $M$ lying in $\mathbb {R}^n$
equipped with a standard Euclidean metric  is $(k,l)$-hyperbolic
if its second quadratic form is everywhere nondegenerate and have
signature $(k,l)$.
\end{Def}

\begin{Def}
We say that a hypersurface $M$ approaches a hypersurface $L$ at
infinity if  for any $\epsilon>0$ there exists an $R>0$ such that
\begin{enumerate}
\item there exists a diffeomorphism $\phi:L\setminus B_R \to M\setminus B_R$,
such that $\|\phi(x)-x\|<\epsilon$ for any $x\in L\setminus
B_{2R}$ and
\item there is a diffeomorphism $\psi$ of the Gauss images of
$L\cap \S_{2R}^{n-1}$ onto the Gauss image of
$\phi(M\cap\S_{2R}^{n-1})\subset L$ such that
$\dist(\psi(x),x)<\epsilon$ in standard metric on $\S^{n-1}$,
where $\S_{2R}^{n-1}$ denotes a sphere of radius $2R$ with center
at the origin.
\end{enumerate}
\end{Def}

The proof of the \thmref{thmmain:1 line exists} goes as follows.
First, using simple topological arguments, we prove that image
under the Gauss map of the hypersurface $M$ does not intersect the
image under the Gauss map of the quadric $\{Q=\e>0\}$ for any
positive $\e$. Second, we note that the interior $U=\{Q>0\}$ of
the cone $K$ can be exhausted by these quadrics. If we suppose
that $M\cap U$ is nonempty and compact, then an analogue of the
Rolle theorem says that there should be a level curve of $Q$ lying
in $U$ and tangent to $M$, namely the level curve $\{Q=\max_{x\in
M}Q(x)\}$. This contradicts to the properties of the Gauss map of
$M$ mentioned above. The assumption of compactness can be dropped
by a slight modification of this arguments.

\subsection{The topological Lemma}\label{ssec:Rolle topology}
The proof of the \thmref{thmmain:1 line exists} starts from a
simple topological lemma. This lemma will be applied further to
the Gauss map of $(k,l)$-hyperbolic hypersurfaces.
\begin{Lem}\label{lem:topol}
Let $M$ be a compact connected manifold with boundary $\partial M$
and $f:M\to N$ be a local diffeomorphism to a compact
simply-connected manifold $N$, $\pi_1(N)=0$. Suppose that the
restriction $f|_{\partial M}$  of $f$ to each connected component
of the boundary  of $M$ is also a diffeomorphism (i.e. also
embedding). Then $f$ is a diffeomorphism of $M$ to $f(M)$.
\end{Lem}
\begin{proof}
An easy case is $\partial M=\emptyset$. In this case $f$ is a
covering, so should be a trivial one (since $\pi_1(N)$, being
trivial, has no nontrivial subgroups).

The general case will be reduced to this case by gluing "hats" to
$M$, thus eliminating the boundary components one-by-one.

Namely, consider a connected component of $\partial M$ (denote it
by $B$). Its image $f(B)$ is a cooriented hypersurface in $N$.
Indeed, $f(B)$ divides any sufficiently small neighborhood of any
its point $f(b)$ into two parts, and one of them is distinguished,
since it is an image of a small neighborhood of $b$ in $M$.
Therefore $N\setminus f(B)$ consists of  two open parts (having
even and odd number of preimages correspondingly). Call the part
which doesn't intersect the image of a sufficiently small
neighborhood of $B$ by "hat". We can glue the "hat" to $M$ along
$B\cong f(B)$: the new manifold will be the union of $B$ and the
"hat" with the neighborhood of $b\in B$ being defined as union of
the connected component of $f^{-1}(U)$ containing $b$ and the
intersection of $U$ and the "hat" (where $U\subset N$ is a small
open ball containing  $f(b)$).

Repeating this operation with all components of the boundary, we
get a new manifold $\widetilde{M}$ without a boundary and
$M\subset  \widetilde{M}$.   The map $f$ extends to
$\widetilde{f}:\widetilde{M}\to N$ in  a trivial manner: we define
it to be an identity on "hats". The map $\widetilde{f}$ satisfies
conditions of the \lemref{lem:topol}, so it is a global
diffeomorphism by the first part of the proof. So $f$ is also a
diffeomorphism, since $f=\widetilde{f}|_M$.
\end{proof}

\begin{Cor}
Let $M$ be a compact connected oriented $(n-1)$-dimensional
submanifold with boundary of $\mathbb{R}^n$, and suppose that its
second fundamental form is everywhere nondegenerate, including the
boundary. Assume that the restriction of the Gauss map of $M$ to
each connected component of the boundary $\partial M$ is
one-to-one. Then the Gauss map of $M$ itself is one-to-one.
\end{Cor}

\subsection{Gauss image of $(k-1,l)$-hyperbolic hypersurface.}
Consider a $(k-1,l)$-hyperbolic  connected closed  hypersurface
$M\subset {\mathbb R}^{k+l+1}$ approaching the quadratic cone
$K=\{Q=0\}$ at infinity, where, as before, $Q(x)=\sum_{i=1}^k
x_i^2-\sum_{j=1}^{l} x_{k+j}^2$ is a quadratic form on $\R^n$,
$n=k+l$. Similarly to the previous Corollary, we prove that the
Gauss image of $M$ coincides with the Gauss image of a $(k-1,l)$-
hyperbolic level hypersurface of $Q$.

We will consider later  the case of a surface $M$ approaching at
infinity the cone $K'=\{(|z|-1)^2=x^2+y^2\}\subset \R^3$ (i.e.
$k-1=l=1$). The \lemref{lem:Gaussimage} is true for this case as
well.

\begin{Lem}\label{lem:Gaussimage}
 Let  $M$ be a  $(k-1,l)$-hyperbolic closed
connected hypersurface approaching $K$. Then the Gauss mapping
maps $M$ diffeomorphically onto the Gauss image of a
$(k-1,l)$-hyperbolic level hypersurface of $Q$.
\end{Lem}

\begin{proof}
The Jacobian of the Gauss map is equal to the Gaussian curvature,
i.e is non-vanishing. Therefore the Gauss map is a local
diffeomorphism. $M$ is not a compact manifold with boundary, so
the \lemref{lem:topol} is not directly applicable.

Consider the compact  $M_R=M\cap B_R$, where $B_R\subset \R^n$ is
a big closed ball of radius $R$ centered at the origin, and denote
by $\partial M_R$ the boundary of $M_R$. Its image $G(\partial
M_R)$ under the Gauss map $C^1$-converge, as $R\to\infty$, to the
Gauss image  of $K\cap \partial B_R$, which coincides with the
Gauss image of the whole cone $K$.
\begin{Lem}\label{lem:exhaustion}
$G(M_R)\cap G(\partial M_{R'})=\emptyset$ for any $R$
and any $R'>R$.
\end{Lem}
For big enough $R'$ this follows immediately from
\lemref{lem:topol} applied to $M_{R'}$ and its Gauss map.  The
\lemref{lem:topol} is applicable since $M_{R'}$ is compact and the
restriction of the Gauss map to its boundary is a diffeomorphism,
due to the condition "approaching at infinity". Therefore this is
true for any $R$.\qed

\begin{Cor} The Gauss image of $M_R$ doesn't intersect the $G(K)$.
\end{Cor}
Indeed, if the intersection is non-empty, then there is a point of
$G(k)$ which lies in the interior of $G(M_{R'})$, where $R'$ is
any number greater than $R$. Therefore this point cannot be a
limit point of $G(\partial M_R)$ as $R\to\infty$. This contradicts
to the condition that $M$ approaches $K$ at infinity.\qed

Therefore the $G(M_R)$ should lie entirely in one of the connected
components into which the $G(K)$ divides the sphere $\S^{n-1}$.
Since $\partial G(M_R)$ converge uniformly to $G(K)$, we conclude
that $G(M)$ is exactly one of them and the Gauss mapping is a
diffeomorphism.

So the last part is to prove that $G(M)$ fall into the right
connected component.  If $k=l$, then all nonsingular level
hypersurfaces of $Q$ are $(k-1,l)$-hyperbolic (since
$(k-1,l)$-hyperbolicity  and $(k,l-1)$-hyperbolicity are the
same), and there is nothing to prove. So we suppose that $k\not
=l$. Then the Gauss image of $K$ divides the sphere $\S^{n-1}$
into two domains of different topological type: the Gauss image of
a $(k-1,l)$-hyperbolic level surface is diffeomorphic to
$\S^{k-1}\times B^l$, and the Gauss image of a
$(k,l-1)$-hyperbolic level hypersurface is diffeomorphic to
$\S^{l-1}\times B^k$. Denote these domains by $D_+$ and $D_-$
correspondingly, so that $D_+=G(\{Q=+1\})$ and $D_-=G(\{Q=-1\})$.
Since the Gauss mapping of $M$ is a diffeomorphism, we get that
$M$ is diffeomorphic to $D_+$ or $D_-$, and our goal is to exclude
the last possibility.

We will prove that $M$ is topologically different from $D_-$. We
will apply the Morse theory to $M$ and the restriction of the
linear functional $f=x_n|_M$ to $M$.

Since $\nabla(x_n)=e_n=(0,0,...,0,1)\in D_-$, the function $f$ has
exactly two critical points on $M$, namely the preimages of $-e_n$
and $e_n$ under the Gauss map of $M$ (which is a diffeomorphism).
The condition of $(k-1,l)$-hyperbolicity means that both critical
points are non-degenerate and  their indices are equal to $k-1$ or
$l$. Let $R$ be a big number and denote by $\widetilde{M}=M\cap
\{x_n<-R\}$. Although the level sets $\{f=c\}=M\cap \{x_n=c\}$ are
not compact, their behavior at infinity is trivial. Indeed, it the
same as of the sections of cone $K$ by the hyperplanes
$\{x_n=c\}$, due to the condition at infinity. Therefore the
standard results of the Morse theory still hold, so the dimension
$h_i(M,\widetilde{M})$ of the group of relative homologies
$H_i(M,\widetilde{M})$ is less than the number of critical points
of index $i$ of the function $f$. The pair $(M,\widetilde{M})$ is
diffeomorphic (through the Gauss maps) to the pair $(Q_{-1},
\widetilde{Q_{-1}})$, where $Q_{-1}=\{Q=-1\}$, and
$\widetilde{Q_{-1}}=Q_{-1}\cap \{x_n<-R\}$, so
$H_i(M,\widetilde{M}=H_i(Q_{-1}, \widetilde{Q_{-1}})$. The latter
can be easily computed to be equal to $1$ for $i=k$ (since
$Q_{-1}/\widetilde{Q_{-1}}\cong\S^{l-1}\vee\S^k$),  and this
contradicts to the fact that $f$ has no critical points of index
$k$ (recall that $k\not =l$).

\end{proof}

\subsection{Rolle Lemma}
Let $Q(x)=x_1^2+...+x_k^2-x_{k+1}^2-...-x_n^2$ be a quadratic form
in $\mathbb{R}^n$. The cone $K=\{Q=0\}$ divides $\mathbb{R}^n$
into two parts, $\{Q>0\}$ and $\{Q<0\}$. We prove in this
paragraph that $M$ does not intersect one of these domains. In
fact the assumption are weaker than before.

\begin{Thm}\label{thm:MoutsideK} Let $M$ be a smooth connected hypersuface $M$ such that
$\dist(x,K)\to 0$ as $M\ni x\to\infty$. Suppose that the Gauss
image  of $M$ is disjoint from the Gauss image of $\{Q=-1\}$. Then
$M$ does not intersect the whole domain $\{Q<0\}$.
\end{Thm}

The proof is a specialization of the following general lemma.
\begin{Lem}\label{lem:maximum}
Let $M\subset {\mathbb R}^n$ be a smooth closed embedded
hypersurface without boundary and suppose that its image under the
Gauss mapping $G_M:M\to {\mathbb S}^{n-1}$ does not intersect a
domain $U\subset{\mathbb S}^{n-1}$. Suppose that $U$ is symmetric
with respect to antipodal map $x\to -x$ of $\S^{n-1}$.

Suppose that on ${\mathbb R}^n$ we are given a function $f$ with
nonnegative only critical values. Suppose that  $\frac{\nabla
f(x)}{\|\nabla f(x)\|}\in U$ if $f(x)<0$ (in other words,
$G(\{f=t\}\subset U$ for any $t<0$). If $f(x)\to 0$ as $M\ni
x\to\infty$, then $f$ is nonnegative on $M$.
\end{Lem}
\begin{proof}
Suppose that $f(x)<0$ for some $x\in M$. Let $x_0\in M$ the  point
of minimum of the restriction of $f$ to $M$. It exists since
$M\cap \{f\le f(M)\}$ is compact and nonempty.  $x_0$ is a
critical point of the restriction of $f$ to $M$. Equivalently,
$T_{x_0}M$ is perpendicular to the nonzero vector $\nabla f(x_0)$.
This means that $\G_S(x_0)\in U$, which is forbidden.
\end{proof}

Apply this Lemma  to the proof of the theorem. A first candidate
for the function $f$ is the $Q$ itself: Gauss images of
$\{Q=t<0\}$ are all equal and do not intersect the Gauss image of
$M$. However, $Q$ itself do not satisfy the conditions of the
\lemref{lem:maximum}: one should slightly adjust $Q$ to ensure
that the restriction of $f$ to $M$ tends to zero at infinity.

Denote   $\sqrt{x_1^2+..+x_k^2}$ by $a$ and
$\sqrt{x_{k+1}^2+...+x_n^2}$ by $b$. Suppose that $Q(x)=a^2-b^2$
take a negative value $-2\e$ at some point $x'\in M$, i.e. that
$M$ intersects the domain $\{Q<-\e\}$. Consider the function
$f_1=\sqrt{a^2+\e}-b$ and denote by $f$ its smoothening: the $f_1$
is not smooth at $b=0$, but one can smoothen $f_1$ without
changing it on $\{f_1<0\}=\{Q<-\e\}$, the only domain interesting
for us. The main point is that the Gauss image of the negative
level hypersurfaces of $f$ is the same as the Gauss image of
$\{Q=-\e\}$. This is most evident in the planar $k=l=1$ case,
where the level curves of $f$ are just translates of $\{Q=-\e\}$.
Indeed, take a point $(x_1,...,x_n)$ such that $t=f(x)<0$. A
simple computation shows that, first, $df(x_1,...,x_n)$ is
proportional to $dQ(x_1,..., x_k,\lambda x_{k+1},...,\lambda x_n)$
for $\lambda=\frac {b+t}{b}$ and, second, $Q(x_1,..., x_k,\lambda
x_{k+1},...,\lambda x_n)=-\e$. This computation proves that the
Gauss images of $\{f=t<0\}$ and $\{Q=-\e\}$ are equal, and also
shows that $f$ has only positive critical values.

One can easily see that different  level curves of $f$ lie on
positive distance one from another, so the only one approaching
$K$ at infinity is the zero level curve $\{f=0\}=\{Q=-\e\}$.
Therefore, since $\dist(x,K)\to 0$ as $M\ni x\to\infty$, the
restriction of $f$ to $M$ tends to zero at infinity.

So $f$ satisfies the conditions of the \lemref{lem:maximum}, and
therefore $M\cap\{Q<-\e\}=M\cap \{f<0\}=\emptyset$, a
contradiction with the choice of $x'$.\qed

\subsection{End of the proof of the \thmref{thmmain:1 line exists}}
The rest of the proof of the \thmref{thmmain:1 line exists} is
just application of the two results proven above, the
\lemref{lem:Gaussimage} and the \thmref{thm:MoutsideK}.

First prove existence of an $l$-dimensional subspace in one of the
domains into which $M$ divides $\mathbb{R}^n$. Suppose first that
$k\not=l$. In this case the quadrics $Q_1=\{Q=1\}$ and
$Q_{-1}=\{Q=-1\}$ have different signatures of the second
quadratic forms: the first one is $(k-1,l)$-hyperbolic, and the
second is $(k,l-1)$-hyperbolic. By \lemref{lem:Gaussimage} the
Gauss image of the $(k-1,l)$-hyperbolic hypersurface $M$ coincide
with the Gauss image of $Q_1$, and is therefore disjoint from the
Gauss image of $Q_{-1}$. So, by \thmref{thm:MoutsideK}, $M$ does
not intersect the domain $\{Q<0\}$, which contains the
$l$-dimensional subspace $\{x_1=...=x_k=0\}$.

If $k=l$, then both $Q_1$ and $Q_{-1}$ are $(k-1,l)$-hyperbolic,
and \lemref{lem:Gaussimage} claims that the Gauss image of $M$
coincides with the Gauss image of one of them. Taking $-Q$ instead
of $Q$ if necessary, we can assume that $G(M)$ coincides with
$G(Q_1)$, and the same arguments hold.

The existence of a $(k-1)$-dimensional affine subspace in the
second part of $\mathbb{R}^n\setminus M$ is evident. By assumption
$M$ approaches the cone $K$ at infinity, so the distance between
$M\setminus B_R$ and $K\setminus B_R$ is less than distance
between $K\setminus B_R$ and  $L=\{x_{k+1}=...=x_n=0\}$ for big
enough ball $B_R$. So take any $(k-1)$-dimensional affine subspace
of $L$ lying outside $B_R$, and it will not intersect $M$.

\section{Projection from the origin}
Starting from this moment we will deal with $(1,1)$-hyperbolic
surfaces in $\mathbb{R}^3$ only. So we will omit the $(1,1)$ and
will call $(1,1)$-hyperbolic surfaces  hyberbolic surfaces.

The \thmref{thmmain:1 line exists} ensures that a hyperbolic
surface approaching the cone $K=\{ x^2+y^2=z^2\}$ at infinity do
not intersect any line passing through the origin and lying in the
domain $\{ x^2+y^2<z^2\}$. We prove here that any line passing
through the origin intersects $M$ at at most two points.
\begin{Thm}\label{thm:central projection}
Let $M\subset {\mathbb R}^3$ be a hyperbolic surface approaching
the standard cone  $K$ at infinity.

Then the restriction to $M$ of the projection $\pi:{\mathbb
R}^3\to{\mathbb S}^2=\{\|x\|=1\}$ is embedding.
\end{Thm}
\subsection{Arnold's formula}
The \thmref{thm:central projection} follows from a remarkable
formula due to Arnold, see \cite{Arnold88}. Consider a generic
smooth hypersurface $M\subset{\mathbb R}P^3$. Denote by
$\#\{M\cap\ell\}$ number of its points of intersections with a
line $\ell$ and by $\sign(M,\ell)$ the number of point $x\in M$
containing the line $\ell$ in their tangent planes counted with
multiplicities. The multiplicity is equal to $"+1"$ if the
Gaussian curvature of $M$ is positive at $x$ and to $"-1"$ if it
is negative at $x$ (if the curvature at $x$ is zero then the
formula for multiplicity is more complicated).
\begin{Lem}[Arnold, 88]\label{lem:Arnold}
For a generic smooth hypersurface $M\subset{\mathbb R}P^3$ the sum
$\#\{M\cap\ell\}+\sign(M,\ell)$ is the same for all $\ell$ and is
equal to the Euler characteristic of $M$
\end{Lem}
\begin{proof}[Sketch of the proof for a semialgebraic $M$ (due to O.Viro)]
 Take out from
$M$ its points of intersection with $\ell$ and compute the Euler
characteristic of the result using Fubini theorem for Euler
characteristic. Namely, the Euler characteristic of $M\setminus
\ell$ is equal to the integral over the space of all planes $L_t$,
$t\in\mathbb{R}P^1$, containing $\ell$ of the Euler characteristic
of $M_t=\{M\setminus\ell\}\cap L_t$. For simplicity, suppose that
each section $M_t$ has at most one singular point (if not, perturb
$\ell$ slightly). Each nonsingular section $M_t$ is a
one-dimensional manifold, so is a union of circles and open
intervals with ends at removed points. Therefore its Euler
characteristic is equal to $-\#\{M\cap\ell\}$ (since Euler
characteristic of a circle is equal to zero). Euler characteristic
of a singular section differs by $+1$ or by $-1$, depending on the
sign of the curvature of $M$ at the singular point. Indeed, if the
curvature at the singular point is negative, then the section has
a self-intersection, so the Euler characteristic drops by $1$. If
the curvature at the singular point is positive, then the section
has an isolated point, and Euler characteristic increases by $1$.

Since the Euler characteristic of $\mathbb{R}P^1$ is equal to
zero, the integration of $-\#\{M\cap\ell\}$ over $\mathbb{R}P^1$
gives zero. So the Euler characteristic of $M\setminus \ell$,
being equal to the integral of the Euler characteristic of $M_t$
over $\mathbb{R}P^1$, is equal to $\sign(M,\ell)$, and the result
follows.\end{proof}

\subsection{Compactification of $M$ and end of the proof of the \thmref{thmmain:2 projection is in}}
We  apply  \lemref{lem:Arnold} to the closure of $M$ in
$\mathbb{R}P^3$. First, we have to show that the closure of $M$ in
${\mathbb R}P^3$ is a smooth surface.
\begin{Lem}\label{lem:projclosure}
The closure $\tilde{M}$ of $M$ in ${\mathbb R}P^3$ is smooth.
\end{Lem}
\begin{proof}
Take  affine coordinates $\tilde{x}=\frac x z,\tilde{y}=\frac y
z,\tilde{w}=\frac 1 z$. We are interested in the points of
$\tilde{M}\cap\{\tilde{w}=0\}$. The first part of the condition
"$M$ approaches $K$ at infinity" implies that $\tilde{M}$
approaches $\{\tilde{x}^2+\tilde{y}^2=1\}$ faster than
$|\tilde{w}|$, so $\tilde{M}$ is smooth at these points. The
second part means that as $x\in M$ tends to $x_0\in
\tilde{M}\cap\{\tilde{w}=0\}$ the limit of tangent planes $T_xM$
exists and is equal to the tangent plane at  $x_0$. This means
$C^1$-smoothness of $\tilde{M}$.
\end{proof}
For the hyperbolic surface $M$ the curvature is always negative.
Therefore  the sign in \lemref{lem:Arnold} is always $"-"$. The
Euler characteristic of $\tilde{M}$ is equal to the Euler
characteristic of $\{x^2+y^2=z^2+w^2\}$ (essentially
\lemref{lem:Gaussimage}), i.e. is equal to zero. So the
\lemref{lem:Arnold} claims in this case that for generic $\ell$
\begin{equation}\label{Arnold for hyperb}
 \#\{M\cap\ell\}=\#\{x\in M\,|\,\ell\subset T_xM\}.
\end{equation}

We want to prove that projection of $M$ to ${\mathbb S}^2$ has no
folds, i.e. is a local diffeomorphism. In other words, we have to
show that no tangent plane  to $M$ passes through the vertex $O$
of the cone. Suppose otherwise and take a plane tangent to $M$ and
passing through  $O$. The normal to this plane lies in the image
of $M$ under the  Gauss map, i.e. in ${\mathbb S}^2\cap\{|z|<\frac
1 {\sqrt(2)}\}$. Equivalently, this plane should intersect the
domain $\{x^2+y^2<z^2\}$. Since this plane contains also the
vertex of the cone, so contains a line $\ell\subset\{x^2+y^2<
z^2\}\cup\{0\}$. By \lemref{thm:MoutsideK} $\ell\cap
\tilde{M}=\emptyset$. Therefore, by compactness of $\tilde{M}$, it
is true for all lines close enough to $\ell$. Moreover, if $\ell$
is contained in a plane $T_xM$, then, due to the nonzero curvature
of $M$ at $x$, any line close enough to $\ell$ is also contained
in some tangent plane to $M$.

So in \ref{Arnold for hyperb} applied to a small perturbation of
$\ell$ the left side is equal to zero, and the right side is at
least one, which is impossible.

The \thmref{thm:central projection} now follows from
\lemref{lem:topol} applied to the restriction of the projection to
$M_R=M\cap B_R$ -- intersection of $M$ with a big enough ball
$B_R$. Indeed, we just proved that the projection is a local
diffeomorphism. Also, the restriction of the projection to the
boundary of $M_R$ is diffeomorphism since the boundary of $M_R$ is
$C^1$-close to $K\cap \{x^2+y^2+z^2=R^2\}$, so is embedded by
projection.\qed

\section{Example}
In this section we provide an example of a hyperbolic closed
connected surface without boundary in $\mathbb{R}^3$ bounding a
domain without lines inside. This surface has an asymptotic
behavior similar to these considered above (namely it approaches
the pseudo-cone $K'$ at infinity), but unlike the surfaces before,
its closure in $\mathbb{R}P^3$ is not smooth and is not hyperbolic
after smoothening.

Construction starts by definition of a affine convex-concave sets.
Consider a hyperbolic surface bounding some domain in
$\mathbb{R}^3$.  At each point it has a direction of positive
sectional curvature and an orthogonal direction of negative
sectional curvature. The affine convex-concave sets come from a
requirement that these directions should not be far from a
vertical (=parallel to $z$-axis) and horizontal (=perpendicular to
$z$-axis) respectively. In other words, we want the horizontal
sectional curvature to be always negative and the vertical
sectional curvature to be always positive. The first requirements
implies that the horizontal sections of a domain bounded by the
hyperbolic surface are convex, and then the second one gives a
concavity-type condition on these sections. We introduce the
affine convex-concave sets as sets satisfying these two
properties, i.e. using only the notion of convexity. This class is
an affine relative of the class of $L$-convex-concave subsets of
$\mathbb{R}P^n$ defined in \cite{NKh2}, and is similarly closed
under surgeries considered there. An analogue of the first part of
the Arnold conjecture can be formulated for convex-concave sets
and for $L$-convex-concave sets (the second part then follows
trivially). We prove it to be true for the first nontrivial case
of $L$-convex-concave sets in \cite{NKh1}, by a rather complicated
considerations.

The first step of our construction in this section is to construct
a counterexample to an analogue of Arnold  conjecture for
convex-concave subsets of $\mathbb{R}^3$. This counterexample is a
so-called strip - a piece of a two-dimensional surface which is at
the same time a convex-concave set. A strip in no way can be
interpreted as a domain bounded by a hyperbolic surface. However,
one can think about it as an interior-less limit of certain
domains bounded by hyperbolic surfaces. The condition of absence
of lines inside a convex-concave body is an open one, so any
convex-concave body close enough to the strip also does not
contains a line inside. The second step consists of a small
perturbation of the set $E$ -  the cone $K'$ with attached strip -
in order to get a convex-concave set which is a genuine domain
bounded by a hyperbolic surface.

This perturbation is based on the fact that the class of
convex-concave sets is closed under taking the  fiberwise affine
linear combinations (by Minkowsky). We replace $E$ by an affine
linear combination of  The result is a genuine domain bounded by a
hyperbolic surface (with an additional property that all its
sections are convex).

\subsection{Affine convex-concave sets}\label{ssec:affine c-c
sets}
We will call by {\em horizontal} everything  in ${\mathbb R}^3$
which is parallel to the coordinate $(x,y)$-plane. For example,
\begin{enumerate}
\item planes $\{z=c\}$ are called horizontal,
\item directions $(a,b,0)$ are called horizontal.
\end{enumerate}

\begin{Def} \label{def:affinecc}
We say that a set $A\subset {\mathbb R}^3$ is convex-concave if
\begin{itemize}
\item its sections by horizontal planes are nonempty, convex and compact and
\item the sections $S_t=A\cap \{z=t\}$ depends in concave way (in Minkowski sense) on
$t$.
\end{itemize}
\end{Def}

The second condition means that for any $t_1<t_2<t_3$ the section
$S_{t_2}$ is contained inside the linear (in Minkowski sense)
combination $\frac{t_2-t_1}{t_3-t_1}S_{t_3}+
\frac{t_3-t_2}{t_3-t_1}S_{t_1}$, i.e. the convex hull of the union
$S_{t_1}\cup S_{t_3}$. It can be reformulated in several possible
ways.

The first equivalent reformulation is that
\begin{itemize}\label{def:affineccv2}
\item[(2')] for any $t_1<t_2<t_3$ any point of the section
$S_{t_2}$ lies on a line intersecting both $S_{t_1}$ and
$S_{t_3}$.
\end{itemize}

Form this follows another reformulation of the same condition.
Namely, it is the requirement that the complement to the
projection of $A$ along any horizontal direction is locally
convex. In other words, if we introduce a coordinate $(w,z)$ on
the plane of projection, then
\begin{itemize}\label{def:affineccv3}
\item[(2'')] the projection should be given by $\pi(A)=\{-\phi_1(z)\le
w\le \phi_2(z)\}$ with both $\phi_1(z)$ and $\phi_2(z)$ being
convex and $-\phi_1(z)\le \phi_2(z)$ for all $z$.
\end{itemize}
Here  "$f$ is a convex function"  means that
$\frac{f(x+a)-f(x)}a-\frac{f(x)-f(x-b)}b\ge 0$ for all $a,b>0$.
For $C^2$-smooth functions this is equivalent to $f''\ge0$, and
for continuous functions it can be defined in more
distribution-like spirit: $f$ is convex if $\int fg''dz\ge 0$ for
any smooth nonnegative function $G$ with compact support or
tending fast enough to zero as $|t|\to \infty$.

\subsubsection{Support function}
Let recall the definition and basic properties of the support
functions of a convex set in $\mathbb{R}^n$ (we will need the case
of $n=2$ only). Let $S\subset\mathbb{R}^n$ be a compact convex
set. Then one can define a support function $F_S(\ell)=\max_{x\in
S}\ell(x)$ on $(\mathbb{R}^n)^*$. This function is clearly
$\mathbb{R}_+$-homogeneous, i.e. $F_S(\alpha\ell)=\alpha
F_S(\ell)$ for any $\alpha>0$ and $\ell\in (\mathbb{R}^n)^*$.
Since $S\subset \{\ell_1(x)\le F_S(\ell_1)\}\cap\{\ell_1(x)\le
F_S(\ell_1)\}$ for any $\ell_1,\ell_2\in (\mathbb{R}^n)^*$, we get
that $F_S(\alpha\ell_1+\beta\ell_2)\le\alpha F_S(\ell_1)+\beta
F_S(\ell_2))$ for any $\alpha,\beta\ge0$. If $S'$ is an another
compact convex subset of $\subset\mathbb{R}^n$, then
$F_S\not=F_{S'}$.

Vice versa,  for  any $F:(\mathbb{R}^2)^*\to \mathbb{R}$
satisfying the previous conditions one can construct a convex
compact figure $S=\cap_{\ell\in(\mathbb{R}^2)^*}\{\ell(x)\le
F(\ell)\}$. One can check that $F_S=F$.

If  we define the Minkowski sum of two sets $A\subset\mathbb{R}^n$
and $B\subset\mathbb{R}^n$ as $A+B=\{a+b|a\in A, b\in B\}$, then
the support function of $A+B$ is equal to the sum of support
functions of $A$ and $B$.

The boundaries of projections in the
\defref{def:affineccv3} are exactly the values of the support
function of $S_t$ on this direction. Namely, if the projection is
defined by $\pi:(x,y,z)\to (ax+by,z)$, then
$\phi_1(z)=\min_{(x,y,z)\in S_z}(-a)x+(-b)y$ and
$\phi_2(z)=\max_{(x,y,z)\in S_z}ax+by$. So the second condition of
the \defref{def:affinecc} means that the support functions of
$S_t$ depend concavely on $t$.

We will further need the case when the boundary of projection
($=$support function) is piecewise

\subsection{Strips}
First  we construct  an unexpected object: a convex-concave set
all horizontal section of which are segments.
\begin{Def}\label{def:strip}
A surface with boundary $\{(x,y,z)\in {\mathbb R}^3\quad |\quad
x=u_1(z)+tf_1(z), y=u_2(z)+tf_2(z), \|t\|\le 1\}$ which is also a
convex-concave set will be called a strip.
\end{Def}

In other words, we parameterize the strip by two curves: one
formed by the middle-points of the segments $M=(u_1(z), u_2(z),
z)$ and another formed by the ends of the segments,
$M_1=(u_1(z)+f_1(z), u_2(z)+f_2(z), z)$

An evident example  of a strip corresponds to linear $f_1$ and
$f_2$ and $u_i\equiv 0$ (so-called {\em degenerate} strip). In
this case the strip is simply a piece of a quadric and contains
two one-parametric families of lines (one of them consists of
horizontal lines). So any point of a degenerate strip lies on a
line intersecting all sections and convex-concavity follows by
\defref{def:affineccv2}.

It turns out that there exist non-degenerate strips, and they
survive some small perturbations, whereas the property to contain
a line doesn't. So we construct a strip not containing a line
inside.

\subsubsection{An unperturbed strip with exactly one line inside}\label{ssec:unpertstrip}

Consider first unperturbed strips, i.e. strips with $u_i(z)\equiv 0$.

We suppose further that $f_i(z)$ are  two linearly independent
solutions of a second order linear differential equation of the
type $y''=g(z)y$.   This is not very restrictive. Indeed, any two
functions are solutions of a differential equation of second order
as soon as their Wronskian is nonzero. But if their Wronskian is
zero at $z=t$, then, assuming some genericity, its sign will
change at this point. It follows (after some computations) that
the projection along direction $(f_1(t),f_2(t),0)$ doesn't satisfy
condition \ref{def:affineccv3}.

Our first step is to build a strip containing exactly one line.
More exact, we will prove that almost all unperturbed strips are
like this.
\begin{Lem}
Suppose that $f_1, f_2$ are linearly independent solutions of a
second order linear differential equation of the type $y''=g(z)y$.
Suppose that the set  $\{x=tf_1(z), y=tf_2(z), \|t\|\le 1\}$
contains another line $\ell$ (apart from the $z$-axis). Then
$g(z)\equiv 0$.
\end{Lem}

\begin{proof}
First, $\ell$ cannot be parallel or intersect  the $z$-axis.
Indeed, in this case our solutions are linearly dependent (since
the equation of a vertical plane is $y=kx$).

So these two lines ($\ell$ and the $z$-axis) are not in the same
plane. After a rotation in $(x,y)$ plane we can assume that
$\ell\in \{x=A\}$, so the line $\ell$ is defined by  equations
$x=A, y=az+b$. Consider the quotient $k(z)=\frac
{f_2(z)}{f_1(z)}$. Then $k'(z)=\frac{W(f_1,f_2)}{f_1^2}=\frac C
{f_1^2}$ since $W(f_1,f_2)\equiv\const$ (the equation $y''=g(z)y$
has no term with $y'$). From the other side, $k(z)=\frac {az+b}A$
is a linear function, so its derivative is a constant. So
$f_1(z)\equiv\const$ and therefore $g(z)\equiv 0$.
\end{proof}

\begin{Lem}\label{lem:strips}
Let $f_i$ be two linearly independent solutions of $y''=g(z)y$.
If,  $g(z)\ge0$ for all $z\in\mathbb R$ then the set
$S=\{x=tf_1(z), y=tf_2(z), \|t\|\le 1\}$ is a strip (i.e. is
convex-concave).
\end{Lem}
\begin{proof}
Projection of the strip in the direction $(-b, a,0)$ is given by
$\pi(A)=\{-|\phi(z)|\le w\le |\phi(z)|\}$, where
$\phi(z)=af_1(z)+bf_2(z)$ is again a solution of the same equation
$y''=g(z)y$. . We have to prove that $|\phi(z)|$ is convex, or,
equivalently, that $\phi(z)''\ge 0$ when $\phi(z)> 0$ and that
$\phi(z)''\le 0$ when $\phi(z)< 0$. But this follows immediately
from the equation and positivity of $g(z)$.
\end{proof}

\subsubsection{Perturbation of a strip}
Take any strip from the \lemref{lem:strips}. We want to perturb
our strip between two levels in such a way that the perturbed set
will be still a strip but will not contain any lines. There is
only one line passing through the unperturbed part and our goal is
to make sure that
\begin{enumerate}
\item the  perturbed part will not
contain this line and
\item that the perturbed strip will remain
convex-concave.
\end{enumerate}

Here is the construction. Take any strip $S$ described in the
\lemref{lem:strips}. Take any $\rho(t)$ such that, first,
$|\rho(z)|\le g(z)$ and, second, ${\rho(z)}\not=\const {g(z)}$.
Take $u_i(z)$ such that $u''_i(z)=\rho(z)f_i(z)$ for $i=1,2$.
Consider the strip $\widetilde{S}=\{(x,y,z)\in {\mathbb R}^3\quad
|\quad x=u_1(z)+tf_1(z), y=u_2(z)+tf_2(z), \|t\|\le 1\}$. In other
words, we shift the segment $S_t$ -- the horizontal section of $S$
-- by vector $(u_1(t), u_2(t),0)$. We prove that the first
condition on $\rho(z)$ implies concave-convexity of
$\widetilde{S}$ and the second implies that $z$-axis$\not\subset
\widetilde{S}$.

\begin{Lem}
$\widetilde{S}$ is  convex-concave.
\end{Lem}
\begin{proof}
The horizontal sections of $S$ are segments, so the first
condition of the \defref{def:affinecc} is satisfied. We check the
second condition in the form of the
\defref{def:affineccv3}. A horizontal projection of $S$ along the
direction $(-b,a,0)$ is given by $\pi(A)=\{\psi(z)-|\phi(z)|\le
w\le \psi(z)+ |\phi(z)|\}$, where $\psi(z)=au_1(z)+bu_2(z)$ and
$\phi(z)=af_1(z)+bf_2(z)$.

We have to check that the boundary of any horizontal projection is
given by convex functions, i.e. that
\begin{itemize}
\item $(\psi(z)+ \phi(z))''\ge 0$ and $(\psi(z)- \phi(z))''\le 0$
when $\phi(z)\ge 0$ and that
\item  $(\psi(z)+ \phi(z))''\le 0$ and $(\psi(z)- \phi(z))''\ge 0$
when $\phi(z)\le 0$
\end{itemize}
In other words, we have to prove  that $(\phi(z)\pm\psi(z))''$ has
the same sign as $\phi(z)$. This is evident since their ratio is
equal to $g\pm\rho$ which is always nonnegative.
\end{proof}

We have to check is that the perturbation is not directed along
the segments, i.e. that $z$-axis does not lie in $\widetilde{S}$.
\begin{Lem}
If $\rho\not =cg$ for a $|c|<1$ then for some $t\in\mathbb R$ the
section $\widetilde{S}\cap\{z=t\}$ of the perturbed strip doesn't
intersect the section $S\cap\{z=t\}$ of the unperturbed one.
\end{Lem}
\begin{proof}
Suppose opposite, i.e. that $u_i=\lambda(z)f_i(z)$. Then $\rho(z)
f_i(z)=u_i''=(\lambda''+\lambda g) f_i+2\lambda'f_i'$.  In other
words, the vector $2\lambda'(f_1',f_2')$ should be proportional to
the vector $(f_1,f_2)$.  This is possible if and only if
$\lambda'\equiv 0$ or the vector $(f_1',f_2')$ is proportional to
the vector $(f_1,f_2)$. The second possibility contradicts to the
linear independence of $f_1$ and $f_2$. The first one means that
$\lambda\equiv\const$, so $\rho$ and $g$ are proportional.
\end{proof}

The perturbation can be made local, i.e. between two levels.
\begin{Lem}\label{lem:strippert}
We can find an even $\rho(z)$ satisfying all previous conditions
and such that $u_i(z)\equiv 0$ for $|z|\ge1/2$.
\end{Lem}
\begin{proof}
Indeed, consider the space $L$ of  even $C^2$-smooth functions
$\rho(z)$ vanishing identically for $|z|\ge1/2$. The functions
$u_i$ solving $u_i''=\rho(z)f_i$ with $\rho\in L$ and initial
conditions $u_i(-1)=u_i'(-1)=0$ are identical zero on $z\le -1/2$
and are linear on $z\ge 1/2$, i.e. $u_i(z)=a_{i1}z+a_{i0}$ for
$z>1/2$. Evidently, $a_{ij}$ depend linearly on $\rho$, so in the
infinite-dimensional space $L$ there is a subspace $L'$, $\codim
L'\le 4$,  of functions corresponding to $a_{ij}=0$.
\end{proof}

\subsubsection{Specification of the strip}
The examples of strips constructed above depend essentially on
$g(z)$ and $\rho(z)$ only  (different choices of $f_i$ and $u_i$
result in strips differing by a linear transformation of
$\mathbb{R}^3$). There is a big degree of freedom in their choice.
Here we impose some additional restrictions on these two functions
and on the choice of $f_i$ and $u_i$, in order to facilitate the
following constructions -- transforming of the strip to a domain
bounded by a hyperbolic surface not containing a line.

\begin{Cor}[of the constructions above]\label{cor:examplestrip}
Take $g(z)$ be an even smooth function identically equal to $0$
for $|z|\ge 1$ and strictly positive otherwise. There is an even
nonzero function  $\rho(z)$ vanishing identically for $|z|\ge 1/2$
and functions $u_1(z),u_2(z),f_1(z),f_2(z)$ such that
\begin{enumerate}
\item $f_i(z)$ are two linearly independent solutions of $f''(z)=g(z)f(z)$ and $u_i''(z)=\rho(z)f_i(z)$;
\item the perturbed strip
$$S=\{(x,y,z)\in {\mathbb R}^3\quad |\quad x=u_1(z)+tf_1(z),\quad
y=u_2(z)+tf_2(z),\quad  \|t\|\le 1\}$$ does not contain lines and
is symmetric  with respect to the rotation $(x,y,z)\to (x,-y,-z)$
of $\mathbb{R}^3$.
\item the part of the strip $S$ lying in $\{|z|\ge 1\}$ is bounded by four rays, and
directions of these rays lie inside the cone $z^2> x^2+y^2$.
\item $f_1^2(\pm2)+f_2^2(\pm2)<1$
\end{enumerate}
\end{Cor}

\begin{proof} Take a $\rho(z)$ as in the \lemref{lem:strippert}.
Take $f_1(z)$ and $f_2(z)$ to be any even and odd correspondingly
solutions of $f''(z)=g(z)f(z)$. Then one can take an even $u_1(z)$
and an odd $u_2(z)$ solutions of equations $u_i''=\rho(z)f_i(z)$.
Since $\rho(z)$ and $g(z)$ are not proportional, the strip $S$
does not contain lines. Together this means that the strip $S$ is
symmetric with respect to the rotation $(x,y,z)\to (x,-y,-z)$ of
$\mathbb{R}^3$.

Since $f_i''(z)=u_i(z)=0$    for $|z|>1$, the boundary of the part
of $S$ lying in $\{|z|>2\}$ is just four rays
$\{(x,y,z)|x=\pm(a_0+a_1|z|), y=\pm b_1|z|, |z|>1\}$. Multiplying
$f_i(z)$ and $u_i(z)$ by a small number (i.e. after a dilatation
of $(x,y)$-plane), we can assume that $a_1^2+b_1^2<1$ and
$f_1^2(\pm2)+f_2^2(\pm2)<1$, as required.
\end{proof}

\subsection{Gluing to the quasi-cone and smoothening}

Here we glue the strip $S$ of the \corref{cor:examplestrip} to the
quasi-cone $K'=\{(x,y,z)\quad|\quad x^2+y^2=(|z|-1)^2, |z|\ge
1\}$. In other words, we construct a convex-concave set $E$ with
horizontal sections  coinciding with sections of $K'$ for
$|z|\ge 2$ and with sections of $S$ for $|z|\le 1$.

Here is how $E$ is constructed. Take the union $E_1$ of $S$ and
$K'$. Horizontal sections of $E_1$  are sometimes segments,
sometimes closed discs and sometimes their unions. Denote by $E$ a
set whose horizontal sections are the convex hulls of the
corresponding horizontal sections of $E_1$. $E$ coincides with
with $S$ for $|z|\le 1$, so in particular doesn't contain a line.
Also, the last two conditions of the \corref{cor:examplestrip}
together guarantee that the part of $S$ lying in $\{|z|\ge 2\}$
lies inside $K'$, i.e. $E$ coincide with $K'$ outside $\{|z|\le
2\}$.

$E$ turns out to be a convex-concave set since its support
function is a maximum of support function of $S$ and a linear
function -- a support function of $K'$ -- overtaking it as
$z\to\infty$, thus convex. Here are the details.
\begin{Lem}
$E$ is a convex-concave set. \end{Lem}
\begin{proof}
All horizontal sections of $E$ are nonempty and convex by
definition. So we have to check that all projections of $E$ are
bounded by graphs of a convex functions, as in
\defref{def:affineccv3}. Let a projection of $E$ be given by
$\pi(E)=\{-\phi_2(z)\le w\le \phi_1(z)\}$. We have to prove that
both $\phi_1(z)$ and $\phi_2(z)$ are convex.

The proof is the same for both $\phi_1(z)$ and $\phi_2(z)$, so we
consider only $\phi_1(z)$. Taking convex hull of sections doesn't
change projection, so $\pi(E)=\pi(E_1)$. Let projections of $S$
and $K'$ be defined by $\pi(S)=\{-\phi^S_2(z)\le w\le
\phi^S_1(z)\}$ and $\pi(K')=\{-\phi^{K'}_2(z)\le w\le
\phi^{K'}_1(z), |z|\ge 1\}$. Then
$\phi_1(z)=\max(\phi^S_1(z),\phi^{K'}_1(z))$ for $|z|\ge 1$ and
$\phi_1(z)=\phi^S_1(z)$ for $|z|<1$.

Let $\widetilde{\phi^{K'}_1}(z)$ be a piece-wise linear function
equal to $\phi^{K'}_1(z)$ for $|z|\ge 1$ and equal to $0$ for
$|z|\le 1$. Trivially  $\widetilde{\phi^{K'}_1}(z)$ is a convex
function.

Note that by choice of $\rho(z)$ in the \corref{cor:examplestrip}
the middle point of the intervals $S\cap \{z=t\}$ lie on the
$z$-axis for $|t|\ge 1/2$, so $\phi^S_1(z)=\phi^S_2(z)\ge 0$ for
$|z|\ge 1/2$. Therefore
$\phi_1(z)=\max(\widetilde{\phi^S_1}(z),\phi^{K'}_1(z))$ for $z\in
[1/2,\infty)$, so is convex on this interval as a maximum of two
convex functions. Similarly $\phi_1(z)$ is convex on
$(-\infty,-1/2]$. By definition $\phi_1(z)$ is a convex function
on $[-1,1]$. Therefore $\phi_1(z)$ is convex on the whole real
line.
\end{proof}

\subsection{Smoothening}

The convex-concave body $E$ built in the previous section doesn't
contain a line but still is not a domain bounded by a hyperbolic
surface. To finish the construction of an example we will smoothen
$E$ and will get a convex-concave domain $D$ bounded by a smooth
hyperbolic surface.

The further constructions are based on the fact that even after
small enough  deformations  $E$ still doesn't contain a line. In
fact, the deformation should be small enough only near the piece
of the strip $S$.
\begin{Lem}\label{lem:smalldeform}
Let $E\subset \mathbb{R}^3$ be some set. Suppose that
$E'\cap\{z=t\}$ is compact for all $t$ and is in an
$\e$-neighborhood of $E\cap\{z=t\}$ for all $-10\le t\le 10$. If
$\e$ is less than some number $c$ depending on $E$, then the set
$E'$ doesn't contain lines.
\end{Lem}
\begin{proof} It is enough to check only non-horizontal lines,
i.e. the lines given by $\ell=\{(a+bz, c+dz,z)|z\in \mathbb{R}\}$.
It is easy to see that the function $\maxdist(\ell,
E)=\max_{x\in\ell\cap\{|z|\le10\}}\dist(x,E)$ achieves its nonzero
minimum $c$ on the set of all non-horizontal lines. Indeed,
$\maxdist(\ell, E)\to\infty$ as $(a,b,c,d)\to\infty$, so there is
a global minimum of $\maxdist(\ell,E)$. Moreover, this minimum is
non-zero since $\maxdist(\ell, E)=0$ would imply that
$\ell\cap\{|z|\le 1\}\subset E\cap\{|z|\le 1\}=S\cap \{|z|\le
1\}$, which is impossible.

If $\e<c$, then $\maxdist(\ell, E')>\maxdist(\ell, E)-\e\ge0$ for
any line $\ell$. This  means that $\ell\not\subset E'\cap\{|z|\le
1\}$, so $\ell\not\subset E'$.
\end{proof}

\subsubsection{Convolution}
The procedure described below are in fact a particular case of a
general method of smoothening of convex-concave sets: though all
considerations are done for the set $E$ constructed in the
previous section, the constructions can be easily generalized for
any convex-concave set with moderate growth of support function as
$|z|\to\infty$. The procedure is a generalization of the
well-known fact that convolution of an integrable function with a
$\mathcal{C}^{\infty}$-smooth function is a
$\mathcal{C}^{\infty}$-smooth function.

There is a well-known operation of taking an affine combination by
Minkowski of two convex sets (by affine combination we mean a
linear combination with positive coefficients sum of which is
equal to $1$). The result of this operation can be described in
two ways:
\begin{enumerate}
\item it is a set of all points which are a convex combination with the same
coefficients of a point in the first set and a point in the second
set.
\item this is a convex set with the support function equal to the
linear combination with the same coefficients of the support
functions of the first and the second sets.
\end{enumerate}

We can apply this operation fiberwise to convex-concave sets.
\begin{Lem}
Let $A$ and $B$ be two convex-concave sets, and let $\lambda_1$
and $\lambda_2$ be two positive numbers, $\lambda_1+\lambda_2=1$.
Then the set $C$ whose horizontal sections $C_z$ are equal to the
$\lambda_1 A_z+\lambda_2 B_z$ is convex-concave.\end{Lem}

Indeed, the sections are convex by definition and support function
of $C_z$, which define boundaries of projections of $C$, are
convex in $z$ as  a sum of two convex in $z$ functions -- the
support functions of $A_z$ and $B_z$.

We will apply a generalization of this operation to the set $E$ of
the previous section. Let take any line $\ell$ perpendicular to
the horizontal planes. Consider a  group of affine transformations
of $\mathbb{R}^3$ generated by translations in vertical directions
and rotations with $\ell$ as an axis. This group $\Gamma$ is just
a cylinder $\Gamma\cong\mathbb{R}\times\mathbb{S}^1$: to
$(z,\phi)$ corresponds a composition $g_{z,\phi}$ of a shift by
$z$ along $\ell$ and rotation by angle $\phi$ around $\ell$. Fix a
standard Lebesgue measure $\mu=\ud z\ud\phi$ on $\Gamma$.

This group acts on the class of convex-concave sets: evidently,
the shifts and rotations around vertical axis of a convex-concave
set gives again a convex-concave set.

As soon as we chose $\ell$, arises a one-to-one correspondence
between convex-concave sets and their support functions
$F_A(z,\ell):\mathbb{R}\times\mathbb{R}^2$: for each fixed $z$
$F_A(z,\ell)$ is just a support function of the section $A_z$,
i.e. $F_A(z,\ell)=\max_{x\in A_z}\ell(x)$. This function is convex
in $z$, is positively homogeneous and is also convex for any fixed
$z$. Vice versa, any such function defines a convex-concave set
(just reconstruct each section separately by its support
function).

Take a $\delta$-like function on $\Gamma$ concentrated near the
identity element of $\Gamma$. More exact, take a function
$\K(z,\phi):\Gamma\to \mathbb{R}$ with a following properties:
\begin{enumerate}
\item $\int_{\Gamma}\K(z,\phi)\ud \mu=1$, and $\int_{|z|,|\phi|\le\e}\K(x,\phi)\ud \mu\ge 1-\e$,
\item $\K(z,\phi)$ is $\mathcal{C}^{\infty}$-smooth and strictly positive,
\item $\K(z,\phi)$ and all its partial derivatives in $z$ decrease exponentially as $|z|\to \infty$
\item $\K(z,\phi)$ is even function of $z$.
\end{enumerate} For example, one can take the
$\K(z,\phi)=C(\e)exp((\cos\phi-z^2)/\e$ with a suitable choice of
the constant $C(\e)$.

We take shifts and rotations of $E$ and define $D_{\e}$ as an
affine combination of the results with weight $\K$. In other
words, we define the convex-concave set $D_{\e}$ as a convolution
of $E$ with $\K(x,\phi)$:
$D=\int_{\Gamma}\K(z,\phi)g_{z,\phi}(E)\ud \mu$. Alternatively,
the support function of the set $D$ is defined by the convolution
$F_D(z,\phi)=\int_{\Gamma}\K(t,\psi)F_E(z-t,\phi-\psi)\ud
t\ud\psi$.

These integrals converge since the support function of $E$ grows
as a linear function of $|z|$: as soon as $|z|\ge 2$, the sections
of $E$ are just circles of radius $|z|-1$.

We claim that
\begin{Thm}
\begin{enumerate}
\item Sections of $D_{\e}$ are strictly convex with nonempty interior.
Moreover, boundaries of sections of $D_{\e}$ are smooth and have
everywhere non-vanishing curvature.
\item $D_{\e}$ is convex-concave domain. Moreover,
boundaries of projections of $D_{\e}$ are smooth and have
everywhere non-vanishing curvature.
\item $D_{\e}$ is close to $E$ in the sense of the
\lemref{lem:smalldeform}. It means that for any $\delta>0$ we can
find  an $\epsilon>0$ such that the sections of domain $D_{\e}$
are in the $\delta$-neighborhoods of the corresponding sections of
$E$. Therefore $D_{\e}$ does not contain lines for a sufficiently
small $\e$.
\item Boundary of $D_{\e}$ is a smooth hyperbolic surface. Moreover, $D_{\e}$ approaches $K'$ at infinity.
\end{enumerate}
\end{Thm}
In other words, the Theorem says that the $D_{\e}$ is an example
we are looking for.

\begin{proof}

Evidently, the support function of $D_{\e}$ is infinitely smooth
as a result of convolution with an infinitely smooth function.
Moreover, since the support function of $E$ convexly and
non-linearly depends on $z$, the support function of $D_{\e}$
will be strictly convex in $z$ and will have everywhere positive
second derivative on $z$.

To prove non-degeneracy  of the curvature of the sections of
$D_{\e}$, one should look on the way to restore a (planar) convex
figure from its support function. Let $F(\ell)$ be a support
function of a planar convex figure $B$, and suppose that $F(\ell)$
is smooth. Then the gradient $\nabla F(\ell)$ is constant on rays
beginning at the origin (i.e. is in fact a mapping $\nabla
F:\S^1\to \mathbb{R}^2$), and its image is exactly the boundary of
$B$. One can easily check that the boundary of $B$ is smooth and
have a nonzero curvature at $\nabla F(\ell_0)$ if the kernel of
the Hessian $H(F)(\ell_0)=\begin{pmatrix}
F_{xx}&F_{xy}\\F_{xy}&F_{yy}\end{pmatrix}(\ell_0)$ is
one-dimensional, i.e. coincides with the line joining the origin
and the point.

We apply this construction to the support functions of sections of
$D_{\e}$. The Hessian of the support function of a section of
$D_{\e}$ is a convolution of $\K$ with the Hessian in coordinates
$(x,y)$ of the support function of $E$. Since the latter is
somewhere nonzero and everywhere positively semi-definite as a
quadratic form, the convolution will be everywhere nonzero and
positive semi-definite. Thus the boundary of sections of $D_{\e}$
are smooth with nonzero curvature.

Moreover, the gradient mapping
$\nabla_{x,y}F_{D_{\e}}(z,\ell):\mathbb{R}\times\S^1\to\partial
D_{\e}$ is a smooth parameterization of the boundary of $D_{\e}$.
Taken together, this means that $D_{\e}$ is convex-concave and its
boundary is smooth with everywhere nonzero curvature.

By standard arguments one can prove that as $\e\to 0$, the result
of a convolution of a function $F_E$ with $\K$ converges to $F_E$
itself. This implies that the set ${D_{\e}}$ lies in a
$\delta$-neighborhood of $E$, as required.

The last claim is that ${D_{\e}}$ approaches $E$ at infinity.
First, the sections $E_z$ of $E$ are circles of radius $|z|-1$ for
$|z|\ge 2$, and therefore $F_E(z,\ell)=(|z|-1)\|\ell\|$ for
$|Z|\ge 2$. Since $\K$ is even as a function of $z$, the
convolution with $\K$ do not changes linear functions of $z$.
Therefore the difference $|F_{D_{\e}}(z,\ell)-F_E(z,\ell)|$
decreases exponentially together with all its derivatives as
$|z|\to\infty$. Therefore the parameterizations of boundaries of
$E$ and ${D_{\e}}$  by the gradient of their support functions as
before are exponentially close, which proves the claim.

\end{proof}
\bibliographystyle{amsplain}

\end{document}